\documentclass[graybox]{svmult}

\usepackage{type1cm}        
\usepackage{makeidx}         
\usepackage{graphicx}        
\usepackage{multicol}        
\usepackage[bottom]{footmisc}

\usepackage{newtxtext}       %
\usepackage{newtxmath}       

\makeindex             

\begin{document}

\title*{Untying The Gordian Knot via Experimental Mathematics}

\author{Yukun Yao and Doron Zeilberger}
\institute{Yukun Yao \at Rutgers University-New Brunswick, 110 Frelinghuysen Rd, Hill Center, Piscataway, NJ 08854, \email{yao@math.rutgers.edu}
\and Doron Zeilberger \at Rutgers University-New Brunswick, 110 Frelinghuysen Rd, Hill Center, Piscataway, NJ 08854, \email{DoronZeil@gmail.com}}
%
%
\maketitle

\abstract{
In this {\it methodological} article on experimental-yet-rigorous enumerative combinatorics,
we use two instructive {\it case studies},
to show that often, just like Alexander the Great before us, the simple, ``cheating'' solution
to a hard problem is the best. So before you spend days (and possibly years) trying to answer
a mathematical question by analyzing and trying to `understand' its structure,  
let your computer generate enough data, and then let it guess the answer.
Often its guess can be  proved by a quick `hand-waving' (yet fully rigorous) `meta-argument'.
Since our purpose is to illustrate a {\it methodology}, we include many details, as well as Maple source-code.}

\section*{Peter Paule}
This article is dedicated to Peter Paule, one of the great pioneers of experimental mathematics
and symbolic computation. In particular, it is greatly inspired by his masterpiece,
co-authored with Manuel Kauers, `The Concrete Tetrahedron' [KP], where a whole chapter
is dedicated to our favorite {\it ansatz}, the $C-$finite ansatz.


\section*{Introduction}

Once upon a time there was a knot that no one could untangle, it was so complicated. Then
came Alexander the Great and, in one second, {\it cut} it with his sword.

Analogously, many mathematical problems are very hard, and the current party line is
that in order for it be considered solved, the solution, or answer, should be given
a logical, rigorous, {\it deductive} proof. 

Suppose that you want to answer the following question: 

\vspace{1mm}\noindent

{\it Find a closed-form formula, as an expression in $n$,
for the real part of the $n$-th complex root of the Riemann zeta function, $\zeta(s)$ .} 

\vspace{1mm}\noindent
Let's call this quantity $a(n)$. Then you compute these real numbers, and find out
that $a(n)=\frac{1}{2}$ for $n \leq 1000$. Later you are told by Andrew Odlyzko that $a(n)=\frac{1}{2}$ for all $1 \leq n \leq 10^{10}$.
Can you conclude that $a(n)=\frac{1}{2}$ for {\it all} $n$? We would, but, at this time of writing, there is
no way to deduce it rigorously, so it remains an open problem. 
It is very possible that one day  it will turn out that $a(n)$ (the real part of the $n$-th complex root of $\zeta(s)$)
belongs to a certain {\it ansatz}, and that checking it for the first $N_0$ cases implies its truth
in general, but this remains to be seen.

There are also frameworks, e.g. {\it Pisot sequences} (see [ESZ], [Z2]), where the {\it inductive } approach fails miserably.

On the other hand, in order to (rigorously) prove that $1^3+2^3+3^3+ \dots + n^3=(n(n+1)/2)^2$, 
for {\it every} positive integer $n$, it suffices to check it
for the five special cases $0 \leq n \leq 4$, since both sides are polynomials of {\bf degree} $4$, hence the
difference is a polynomial of degree $\leq 4$, given by five `degrees of freedom'. 

This is an example of what is called  the  `$N_0$ principle'. In the case of a polynomial identity (like this one),
$N_0$ is simply the degree plus one.

But our favorite {\it ansatz} is the $C$-finite ansatz. A sequence of numbers $\{a(n)\}$ ($0 \leq n < \infty$)
is $C$-finite if it satisfies a {\it linear recurrence equation with constant coefficients}.
For example the Fibonacci sequence that satisfies $F(n)-F(n-1)-F(n-2)=0$ for $n \geq 2$.

The $C$-finite ansatz is beautifully described in Chapter 4 of  the masterpiece
`The Concrete Tetrahedron' ([KP]), by Manuel Kauers and Peter Paule, and discussed at length in [Z1].

Here the  `$N_0$ principle' also holds (see [Z3]), i.e. by looking at the `big picture' one can determine
{\it a priori}, a positive integer, often not that large, such that checking that $a(n)=b(n)$ for $1 \leq n \leq N_0$
implies that $a(n)=b(n)$ for all $n>0$.

A sequence $\{a(n)\}_{n=0}^{\infty}$ is $C$-finite if and only if its (ordinary) {\it generating function}
$f(t):=\sum_{n=0}^{\infty} a(n)\,t^n$ is a {\bf rational function} of $t$, i.e.
$f(t)=P(t)/Q(t)$ for some {\it polynomials} $P(t)$ and $Q(t)$. 
For example, famously, the generating function of the Fibonacci
sequence is $t/(1-t-t^2)$.

Phrased in terms of generating functions, the $C$-finite ansatz is the subject of chapter 4 of yet another
masterpiece, Richard Stanley's `Enumerative Combinatorics' (volume 1) ([S]). There it is shown, using the
`transfer matrix method' (that originated in physics), that in many combinatorial situations, where
there are finitely many states, one is guaranteed, {\it a priori}, that the generating function is rational.

Alas, finding this transfer matrix, at each specific case, is not easy!
The human has to first figure out the set of states, and then using human ingenuity, figure out
how they interact.

A better way is to automate it. Let the computer do the research, and using `symbolic dynamical programming',
the computer, automatically, finds the set of states, and constructs, {\it all by itself} (without any human
pre-processing) the set of states and the transfer matrix. But this may not be so efficient for two reasons.
First, at the very end, one has to invert a matrix with {\it symbolic} entries, hence compute 
symbolic determinants, that is time-consuming. Second, setting up the `infra-structure' and writing a program
that would enable the computer to do `machine-learning' can be very daunting.

In this article, we will describe two {\it case studies} where, by `general nonsense', we know
that the generating functions are rational, and it is easy to bound the degree of the denominator
(alias the order of the recurrence satisfied by the sequence).  
Hence a simple-minded, {\it empirical}, approach of computing the first few terms and then
`fitting' a recurrence (equivalently rational function) is possible.

The first case-study concerns
counting spanning trees in families of grid-graphs, studied by Paul Raff ([R]), and F.J. Faase ([F]).
In their research, the human first analyzes the intricate combinatorics, manually sets up
the transfer matrix, and only at the end lets a computer-algebra system evaluate the symbolic determinant.

Our key observation, that enabled us to `cut the Gordian knot' is that  the terms of the studied sequences are
expressible as {\it numerical} determinants. Since computing numerical determinants is so fast,
it is easy to compute sufficiently many terms, and then fit the data into a rational function.
Since we easily have an upper bound for the degree of the denominator of the rational function,
everything is rigorous.

The second case-study is computing generating functions for sequences of determinants of `almost diagonal Toeplitz matrices'.
Here, in addition to the `naive' approach of cranking enough data and then fitting it into a rational function, we
also describe the `symbolic dynamical programming method', that surprisingly is faster for the range
of examples that we considered. But we believe that for sufficiently large cases, the naive approach will
eventually be more efficient, since the `deductive' approach works equally well for the analogous problem of
finding the sequence of permanents of these almost diagonal Toeplitz matrices, for which the naive approach will soon be
intractable.

This article may be viewed as a {\it tutorial}, hence we include lots of implementation details, and
Maple code. We hope that it will  inspire readers (and their computers!) to apply it in other situations


\section*{Accompanying Maple Packages}

This article is accompanied by three Maple packages,
{\tt GFMatrix.txt}, {\tt JointConductance.txt}, and {\tt SpanningTrees.txt}, all available from
the url

\bigskip

{\tt http://sites.math.rutgers.edu/\~{}zeilberg/mamarim/mamarimhtml/gordian.html}

\bigskip

In that page there are also links to numerous sample input and output files.

\section*{The human approach to enumerating spanning trees of grid graphs}

In order to illustrate the advantage of ``keeping it simple'', we will review  the human approach to 
the enumeration task that we will later redo using the `Gordian knot' way.
While the human approach is definitely interesting for its own sake, it is rather painful.

Our goal is to enumerate the number of spanning trees 
in certain families of graphs, 
notably grid graphs and their generalizations. Let's examine  Paul Raff's interesting approach
described in his paper {\it Spanning Trees in Grid Graph} [R]. Raff's approach was inspired by 
the pioneering  work of F. J. Faase  ([F]).

The goal is to find generating functions that
enumerate spanning trees in grid graphs and the product of an arbitrary graph and a path or a cycle. 

Grid graphs have two parameters, let's call them $k$ and $n$. 
For a $k \times n$ grid graph, let's think of $k$ as {\it fixed} while $n$ is the discrete input variable of interest.

{\bf Definition} The $k \times n$ grid graph $G_k(n)$ is the following graph given in terms of its vertex
set $V$ and edge set $E$:
$$
V = \{v_{ij}|1 \leq i \leq k, 1 \leq j \leq n\},
$$
$$
S = \{\{v_{ij}, v_{i'j'}\}| |i-i'|+|j-j'|=1\}.
$$
The main idea in the human approach is to consider  the collection of set-partitions  of
$[k] = \{1,2,\dots,k\}$ and figure out  the transition when we extend a $k \times n$ grid graph to a $k \times (n+1)$ one.

Let $\mathcal{B}_k$ be the collection  of all set-partitions of $[k]$.  
$B_k = |\mathcal{B}_k|$ is called the $k$-th Bell number. Famously, the exponential generating function
of $B_k$, namely $\sum_{k=0}^{\infty} \frac{B_k}{k!}\, t^k$, equals $e^{e^t-1}$.

A lexicographic ordering on $\mathcal{B}_k$ is defined as follows:

{\bf Definition} Given two partitions $P_1$ and $P_2$ of $[k]$, for $i \in [k]$, let $X_i$ be the block of $P_1$ containing $i$ and $Y_i$ be the 
block of $P_2$ containing $i$. Let $j$ be the minimum number such that $X_i \neq Y_i$. Then $P_1 < P_2$ iff

1. $|P_1| < |P_2|$ or

2. $|P_1| = |P_2|$ and $X_j \prec Y_j$ where $\prec$ denotes the normal lexicographic order. 

For example, here is the ordering for $k=3$:
$$
\mathcal{B}_3 = \{\{\{1,2,3\}\}, \{\{1\}, \{2,3\}, \{\{1,2\}, \{3\}\}, \{\{1,3\}, \{2\}\}, \{\{1\}, \{2\}, \{3\}\}\} \quad .
$$
For simplicity, we can rewrite it as follows:
$$
\mathcal{B}_3 = \{123, 1/23, 12/3, 13/2, 1/2/3\}.
$$
{\bf Definition} Given a spanning forest $F$ of $G_k(n)$, the partition induced by $F$ is obtained from the equivalence relation

\centerline{$i \sim j \Longleftrightarrow v_{n,i}, v_{n,j} $ are in the same component of $F$.}

For example, the partition induced by any spanning tree of $G_k(n)$ is $123\dots k$ because by definition,
in a spanning tree, all $v_{n,i}, 1 \leq i \leq k$ are in the same component.
For the other extreme,
where every component only consists of one vertex, the corresponding set-partition is
$1/2/3/\dots /k-1/k$ because no two $v_{n,i}, v_{n,j}$ are in the same component for $1 \leq i<j \leq k$.

{\bf Definition} Given a spanning forest $F$ of $G_k(n)$ and a set-partition $P$ of $[k]$, we say that $F$ is consistent with $P$ if:

1. The number of trees in $F$ is precisely $|P|$.

2. $P$ is the partition induced by $F$.

Let $E_n$ be the set of edges $E(G_k(n)) \backslash E(G_k(n-1))$, then $E_n$ has $2k-1$ members.

Given a forest $F$ of $G_k(n-1)$ and some subset $X \subseteq E_n$, we can combine them to 
get a forest of $G_k(n)$ as follows. 
We just need to know  how many subsets of $E_n$ can transfer a forest consistent with some partition to a forest consistent with another partition. 
This leads to the following definition:

{\bf Definition} Given two partitions $P_1$ and $P_2$ in $\mathcal{B}_k$, a subset $X \subseteq E_n$ transfers from $P_1$ to $P_2$ if a forest 
consistent with $P_1$ becomes a forest consistent with $P_2$ after the addition of $X$. In this case, we write $X \diamond P_1 = P_2$.

With the above definitions, it is natural to define a $B_k \times B_k$ transfer matrix $A_k$ by the following:
$$
A_k(i,j) = | \{A \subseteq E_{n+1} | A \diamond P_j = P_i \} |.
$$
Let's look at the $k=2$ case as an example. We have 
$$
\mathcal{B}_2 = \{12, 1/2\}, E_{n+1} = \{\{v_{1,n}, v_{1,n+1}\}, \{v_{2,n}, v_{2,n+1}\}, \{v_{1,n+1}, v_{2,n+1}\}\}.
$$
For simplicity, let's call the edges in $E_{n+1}$ $e_1, e_2, e_3$. Then to transfer the set-partition $P_1 = 12$ to itself, 
we have the following three ways: $\{e_1, e_2\}, \{e_1, e_3\}, \{e_2, e_3\}$. 
In order to transfer the partition $P_2=1/2$ into $P_1$, we  only have one way, namely:
$\{e_1, e_2, e_3\}$. Similarly, there are two ways to transfer $P_1$ to $P_2$ and one way to transfer $P_2$ to itself 
Hence the transfer matrix is the following $2 \times 2$ matrix:
$$
A =
\begin{bmatrix}
 3 & 1 \\
 2 & 1
\end{bmatrix}.
$$
Let $T_1(n), T_2(n)$ be the number of forests of $G_k(n)$ which are consistent with the partitions 
$P_1$ and  $P_2$, respectively. Let
$$
v_n =
\begin{bmatrix}
 T_1(n) \\
 T_2(n)
\end{bmatrix} \quad ,
$$
then
$$
 v_n =Av_{n-1} \quad .
$$
The characteristic polynomial of $A$ is
$$
\chi_\lambda(A) = \lambda^2-4\lambda+1.
$$
By the Cayley-Hamilton Theorem, $A$ satisfies 
$$
A^2-4A+1=0.
$$
Hence the recurrence relation for $T_1(n)$ is
$$
T_1(n) = 4T_1(n-1) - T_1(n-2),
$$
the sequence is $\{1, 4, 15, 56, 209, 780, 2911, 10864, 40545, 151316, \dots \}$ (OEIS A001353) and the generating function is 
$$
\frac{x}{1-4x+x^2}.
$$
Similarly, for the $k=3$ case, the transfer matrix
$$
A_3 =
\begin{bmatrix}
 8 & 3 & 3 & 4 & 1 \\
 4 & 3 & 2 & 2 & 1 \\
 4 & 2 & 3 & 2 & 1 \\
 1 & 0 & 0 & 1 & 0 \\
 3 & 2 & 2 & 2 & 1 
\end{bmatrix}.
$$
The transfer matrix method can be generalized to general graphs of the form $G \times P_n$, especially cylinder graphs. 

As one can see, 
we had to think very hard.
First we had to establish a `canonical'
ordering over set-partitions, then
define the consistence between partitions and forests, then look for the transfer matrix and finally worry about initial conditions. 

Rather than think so hard, let's compute sufficiently many terms of the enumeration sequence, and try to guess a linear
recurrence equation with constant coefficients, that would be provable {\it a posteriori} just because we know that {\it there exists}  a transfer matrix
without worrying about finding it explicitly. But how do we generate sufficiently many terms?
Luckily, we can use the celebrated {\bf Matrix Tree Theorem}. 

\subsection*{The  Matrix Tree Theorem}

{\bf Matrix Tree Theorem} If $A = (a_{ij})$ is the adjacency matrix of an arbitrary graph $G$, then the number of spanning trees is equal to the 
determinant of any co-factor of the Laplacian matrix $L$ of $G$, where 

$$ L = 
\begin{bmatrix}
    a_{12}+\dots+a_{1n} & -a_{12} &  \dots  & -a_{1,n} \\
    -a_{21} & a_{21}+\dots+a_{2n}  & \dots & -a_{2,n} \\
    \vdots & \vdots  & \ddots & \vdots \\
    -a_{n1} & -a_{n2} & \dots  & a_{n1}+\dots+a_{n,n-1}
\end{bmatrix}.
$$

For instance, taking the $(n,n)$ co-factor, we have that the number of spanning trees of $G$ equals
$$
\begin{vmatrix}
    a_{12}+\dots+a_{1n} & -a_{12} &  \dots  & -a_{1,n-1} \\
    -a_{21} & a_{21}+\dots+a_{2n}  & \dots & -a_{2,n-1} \\
    \vdots & \vdots  & \ddots & \vdots \\
    -a_{n-1,1} & -a_{n-1,2} & \dots  & a_{n-1,1}+\dots+a_{n-1,n}
\end{vmatrix}.
$$
Since computing determinants for numeric matrices is very fast, we can find the generating functions for the number of spanning trees in grid graphs and more generalized graphs by experimental methods, using the C-finite ansatz.

\section*{The GuessRec Maple procedure}

Our engine is the Maple procedure  {\tt GuessRec(L)} that resides 
in the Maple packages accompanying this article.
We used the `vanilla', straightforward, linear algebra  approach for guessing, using {\it undetermined coefficients}.
A more efficient way is via the celebrated Berlekamp-Massey algorithm ([Wi]). Since the guessing part is not the {\it bottle-neck} of
our approach  ( it is rather the data-generation part), we preferred to keep it simple.

Naturally, we need to collect enough data. The input is the data (given as a list) 
and the output is a conjectured recurrence relation derived from that data.

Procedure {\tt GuessRec(L)}  inputs a list, {\tt L}, and attempts to output a linear recurrence equation with constant coefficients
satisfied by the list. It is based on procedure  {\tt GuessRec1(L,d)} that looks for such a recurrence of order $d$.

The output of   {\tt GuessRec1(L,d)} consists of the
the list of initial $d$ values (`initial conditions')
and the recurrence equation represented as a list. 
For instance, if the input is $L = [1,1,1,1,1,1]$ and $d=1$, then the output will be $[[1],[1]]$; 
if the input is $L=[1, 4, 15, 56, 209, 780, 2911, 10864, 40545, 151316]$ as the $k=2$ case for grid graphs and $d=2$, then the output will be 
$[[1, 4], [4, -1]]$. This means that our sequence satisfies the recurrence $a(n)=4a(n-1)-a(n-2)$, subject to the initial conditions
$a(0)=1,a(1)=4$.

Here is the Maple code:

\bigskip

\vspace{1mm}\noindent

{\obeylines
{\tt
GuessRec1:=proc(L,d) local eq,var,a,i,n:
if nops(L)<=2*d+2 then
\quad    print(`The list must be of size >=`, 2*d+3 ):
\quad    RETURN(FAIL):
fi:
var:=$\{$seq(a[i],i=1..d)$\}$:
eq:=$\{$seq(L[n]-add(a[i]*L[n-i],i=1..d),n=d+1..nops(L))$\}$:
var:=solve(eq,var):
if var=NULL then
\quad    RETURN(FAIL):
else
\quad   RETURN([[op(1..d,L)],[seq(subs(var,a[i]),i=1..d)]]):
fi:
end:
}
}

\bigskip

The  idea is that having a long enough list $L$ $(|L|>2d+2)$  of data, we use the data after the $d$-th one to 
discover whether there exists a linear recurrence relation, the first $d$ data points being the initial condition. 
With the unknowns $a_1, a_2, \dots, a_d $, we have a linear systems of no less than $d+3$ equations. If there is a solution, 
it is extremely likely that the recurrence relation holds in general. 
The first list of length $d$ in the output constitutes the list of initial conditions while the second list, $R$,
codes the linear recurrence, where $[R[1], \dots R[d]]$ stands for the following recurrence:
$$
L[n] = \sum_{i=1}^d R[i]L[n-i].
$$

Here is the Maple procedure {\tt GuessRec(L)}:

\bigskip

{\obeylines
{\tt
GuessRec:=proc(L) local gu,d:
for d from 1 to trunc(nops(L)/2)-2 do
\quad  gu:=GuessRec1(L,d):
\quad  if gu<>FAIL then
  \quad \quad  RETURN(gu):
       fi:
od:
FAIL:
end:
}
}

\bigskip

This procedure inputs a sequence $L$ and tries to guess a recurrence equation with constant coefficients satisfying it. 
It returns the initial values and the recurrence equation as a pair of lists. Since the length of $L$ is limited, the maximum degree of the recurrence 
cannot be more than $\lfloor |L|/2-2 \rfloor$. With this procedure, we just need to input 
$L=[1, 4, 15, 56, 209, 780, 2911, 10864, 40545, 151316]$  to get the recurrence (and initial conditions) $[[1, 4], [4, -1]]$.

Once the recurrence relation, let's call it {\tt S},
is discovered, procedure {\tt CtoR(S,t)} 
finds the generating function for the sequence. 
Here is the Maple code:

\bigskip

{\obeylines
{\tt
CtoR:=proc(S,t) local D1,i,N1,L1,f,f1,L:
if not (type(S,list) and  nops(S)=2 and type(S[1],list) and type(S[2],list)
\quad        and nops(S[1])=nops(S[2]) and type(t, symbol) ) then
\quad  print(`Bad input`):
\quad   RETURN(FAIL):
fi:
D1:=1-add(S[2][i]*t**i,i=1..nops(S[2])):
N1:=add(S[1][i]*t**(i-1),i=1..nops(S[1])):
L1:=expand(D1*N1):
L1:=add(coeff(L1,t,i)*t**i,i=0..nops(S[1])-1):
f:=L1/D1:
L:=degree(D1,t)+10:
f1:=taylor(f,t=0,L+1):
if expand([seq(coeff(f1,t,i),i=0..L)])<>expand(SeqFromRec(S,L+1)) then
print([seq(coeff(f1,t,i),i=0..L)],SeqFromRec(S,L+1)):
\quad  RETURN(FAIL):
else
\quad  RETURN(f):
fi:
end:
}
}

\vspace{1mm}\noindent

Procedure {\tt SeqFromRec} used above (see the package) simply generates many terms using the recurrence.

Procedure {\tt CtoR(S,t)} outputs the rational function in $t$, 
whose coefficients are the members of the C-finite sequence $S$. For example: 
$$
{\tt CtoR([[1,1],[1,1]],t)} = \frac{1}{-t^2-t+1}.
$$
Briefly, the idea is that the denominator of the rational function can be easily determined by the recurrence relation and we use the initial 
condition to find the starting terms of the generating function,
then multiply it by the denominator, yielding the numerator.

\section*{Application of GuessRec for enumerating  spanning trees of grid graphs and  $G \times P_n$}

With the powerful  procedures {\tt GuessRec} and {\tt CtoR}, we are able to find generating functions 
for the number of spanning trees of generalized graphs of the form $G \times P_n$. We will illustrate 
the application of {\tt GuessRec} to finding the generating function for the number of spanning trees in grid graphs. 

First, using procedure {\tt GridMN(k,n)}, we get the  $k \times n$ grid graph.

Then, procedure {\tt SpFn} uses the Matrix Tree Theorem
to evaluate the determinant of the co-factor of the Laplacian matrix of the grid graph which is the number of spanning trees in this particular graph. 
For a fixed $k$, we need to generate a sufficiently long list of data for the number of spanning trees in 
$G_k(n), n \in [l(k), u(k)]$. 
The lower bound $l(k)$ can't be too small since the first several terms are the initial condition; the upper bound $u(k)$ can't be too small as well 
since we need sufficient data to obtain the recurrence relation. 
Notice that there is a symmetry for the recurrence relation, and
to take advantage of this fact,  we modified {\tt GuessRec} to get the more efficient {\tt GuessSymRec} 
(requiring less data).
Once the recurrence relation, and the initial conditions, are given, 
applying {\tt CtoR(S,t)} will give the desirable generating function, that, of course, is a rational function of $t$.
All the above is incorporated in 
procedure {\tt GFGridKN(k,t)} which inputs a positive integer $k$ and a symbol $t$, 
and outputs the generating function whose coefficient of $t^n$ is the number of spanning trees in $G_k(n)$, i.e. if we let 
$s(k,n)$ be the number of spanning trees in $G_k(n)$, the generating function 
$$
F_k(t) = \sum_{n=0}^{\infty} s(k,n) t^n.
$$
We now list  the generating functions $F_k(t)$ for $1 \leq k \leq 7$:
Except for $k=7$, these were already found by Raff[R] and Faase[F], but it is
reassuring that, using our new approach, we got the same output. The case $k=7$ seems to be new.

\bigskip

{\bf Theorem 1} The generating function for the number of spanning trees in $G_1(n)$ is:
$$
F_1(t) = \frac {t}{1-t}.
$$

\bigskip

{\bf Theorem 2} The generating function for the number of spanning trees in $G_2(n)$ is:
$$
F_2 = \frac {t}{{t}^{2}-4\,t+1}.
$$

\bigskip

{\bf Theorem 3} The generating function for the number of spanning trees in $G_3(n)$ is:
$$
F_3 = \frac {-{t}^{3}+t}{{t}^{4}-15\,{t}^{3}+32\,{t}^{2}-15\,t+1}.
$$

\bigskip

{\bf Theorem 4} The generating function for the number of spanning trees in $G_4(n)$ is:
$$
F_4 = \frac {{t}^{7}-49\,{t}^{5}+112\,{t}^{4}-49\,{t}^{3}+t}{{t}^{8}-56\,{t
}^{7}+672\,{t}^{6}-2632\,{t}^{5}+4094\,{t}^{4}-2632\,{t}^{3}+672\,{t}^
{2}-56\,t+1}.
$$
For $5 \leq k \leq 7$, since the formulas are too long, we present their numerators and denominators separately.

\bigskip

{\bf Theorem 5} The generating function for the number of spanning trees in $G_5(n)$ is:
$$
F_5 = \frac{N_5}{D_5}
$$
where
$$
N_5 = -{t}^{15}+1440\,{t}^{13}-26752\,{t}^{12}+185889\,{t}^{11}-574750\,{t}^
{10}+708928\,{t}^{9}-708928\,{t}^{7} 
$$
$$
+574750\,{t}^{6}-185889\,{t}^{5}+
26752\,{t}^{4}-1440\,{t}^{3}+t,
$$

$$
D_5 = {t}^{16}-209\,{t}^{15}+11936\,{t}^{14}-274208\,{t}^{13}+3112032\,{t}^{
12}-19456019\,{t}^{11}+70651107\,{t}^{10}-152325888\,{t}^{9}
$$
$$
+196664896\,{t}^{8}-152325888\,{t}^{7}+70651107\,{t}^{6}-19456019\,{t}^{5}+
3112032\,{t}^{4}-274208\,{t}^{3}+11936\,{t}^{2}-209\,t+1.
$$

\bigskip

{\bf Theorem 6} The generating function for the number of spanning trees in $G_6(n)$ is:
$$
F_6 = \frac{N_6}{D_6}
$$
where
$$
N_6 = {t}^{31}-33359\,{t}^{29}+3642600\,{t}^{28}-173371343\,{t}^{27}+
4540320720\,{t}^{26}-70164186331\,{t}^{25}
$$
$$
+634164906960\,{t}^{24}-
2844883304348\,{t}^{23}-1842793012320\,{t}^{22}+104844096982372\,{t}^{
21}
$$
$$
-678752492380560\,{t}^{20}+2471590551535210\,{t}^{19}-
5926092273213840\,{t}^{18}+9869538714631398\,{t}^{17}
$$
$$
-11674018886109840\,{t}^{16}+9869538714631398\,{t}^{15}-
5926092273213840\,{t}^{14}+2471590551535210\,{t}^{13}
$$
$$
-678752492380560
\,{t}^{12}+104844096982372\,{t}^{11}-1842793012320\,{t}^{10}-
2844883304348\,{t}^{9}
$$
$$
+634164906960\,{t}^{8}-70164186331\,{t}^{7}+
4540320720\,{t}^{6}-173371343\,{t}^{5}+3642600\,{t}^{4}-33359\,{t}^{3}
+t,
$$

$$
D_6 = {t}^{32}-780\,{t}^{31}+194881\,{t}^{30}-22377420\,{t}^{29}+1419219792
\,{t}^{28}-55284715980\,{t}^{27}+1410775106597\,{t}^{26}
$$
$$
-24574215822780\,{t}^{25}+300429297446885\,{t}^{24}-2629946465331120\,{
t}^{23}+16741727755133760\,{t}^{22}
$$
$$
-78475174345180080\,{t}^{21}+
273689714665707178\,{t}^{20}-716370537293731320\,{t}^{19}
$$
$$
+1417056251105102122\,{t}^{18}-2129255507292156360\,{t}^{17}+
2437932520099475424\,{t}^{16}
$$
$$
-2129255507292156360\,{t}^{15}+
1417056251105102122\,{t}^{14}-716370537293731320\,{t}^{13}
$$
$$
+273689714665707178\,{t}^{12}-78475174345180080\,{t}^{11}+
16741727755133760\,{t}^{10}-2629946465331120\,{t}^{9}
$$
$$
+300429297446885
\,{t}^{8}-24574215822780\,{t}^{7}+1410775106597\,{t}^{6}-55284715980\,
{t}^{5}
$$
$$
+1419219792\,{t}^{4}-22377420\,{t}^{3}+194881\,{t}^{2}-780\,t+1.
$$
{\bf Theorem 7} The generating function for the number of spanning trees in $G_7(n)$ is:
$$
F_7 = \frac{N_7}{D_7}
$$
where
$$
N_7 = -{t}^{47}-142\,{t}^{46}+661245\,{t}^{45}-279917500\,{t}^{44}+
53184503243\,{t}^{43}-5570891154842\,{t}^{42}
$$
$$
+341638600598298\,{t}^{41
}-11886702497030032\,{t}^{40}+164458937576610742\,{t}^{39}
$$
$$
+4371158470492451828\,{t}^{38}-288737344956855301342\,{t}^{37}+
7736513993329973661368\,{t}^{36}
$$
$$
-131582338768322853956994\,{t}^{35}+
1573202877300834187134466\,{t}^{34}
$$
$$
-13805721749199518460916737\,{t}^{
33}+90975567796174070740787232\,{t}^{32}
$$
$$
-455915282590547643587452175\,
{t}^{31}+1747901867578637315747826286\,{t}^{30}
$$
$$
-5126323837327170557921412877\,{t}^{29}+11416779122947828869806142972\,
{t}^{28}
$$
$$
-18924703166237080216745900796\,{t}^{27}+
22194247945745188489023284104\,{t}^{26}
$$
$$
-15563815847174688069871470516
\,{t}^{25}+15563815847174688069871470516\,{t}^{23}
$$
$$
-22194247945745188489023284104\,{t}^{22}+18924703166237080216745900796
\,{t}^{21}
$$
$$
-11416779122947828869806142972\,{t}^{20}+
5126323837327170557921412877\,{t}^{19}
$$
$$
-1747901867578637315747826286\,{
t}^{18}+455915282590547643587452175\,{t}^{17}
$$
$$
-90975567796174070740787232\,{t}^{16}+13805721749199518460916737\,{t}^{
15}
$$
$$
-1573202877300834187134466\,{t}^{14}+131582338768322853956994\,{t}^
{13}-7736513993329973661368\,{t}^{12}
$$
$$
+288737344956855301342\,{t}^{11}-
4371158470492451828\,{t}^{10}-164458937576610742\,{t}^{9}
$$
$$+11886702497030032\,{t}^{8}-341638600598298\,{t}^{7}+5570891154842\,{t}
^{6}-53184503243\,{t}^{5}
$$
$$
+279917500\,{t}^{4}-661245\,{t}^{3}+142\,{t}^
{2}+t,
$$

$$
D_7 = {t}^{48}-2769\,{t}^{47}+2630641\,{t}^{46}-1195782497\,{t}^{45}+
305993127089\,{t}^{44}-48551559344145\,{t}^{43}
$$
$$
+5083730101530753\,{t}^
{42}-366971376492201338\,{t}^{41}+18871718211768417242\,{t}^{40}
$$
$$
-709234610141846974874\,{t}^{39}+19874722637854592209338\,{t}^{38}-
422023241997789381263002\,{t}^{37}
$$
$$
+6880098547452856483997402\,{t}^{36}
-87057778313447181201990522\,{t}^{35}
$$
$$
+862879164715733847737203343\,{t}
^{34}-6750900711491569851736413311\,{t}^{33}
$$
$$
+41958615314622858303912597215\,{t}^{32}-208258356862493902206466194607
\,{t}^{31}
$$
$$
+828959040281722890327985220255\,{t}^{30}-
2654944041424536277948746010303\,{t}^{29}
$$
$$
+6859440538554030239641036025103\,{t}^{28}-
14324708604336971207868317957868\,{t}^{27}
$$
$$
+24214587194571650834572683444012\,{t}^{26}-
33166490975387358866518005011884\,{t}^{25}
$$
$$
+36830850383375837481096026357868\,{t}^{24}-
33166490975387358866518005011884\,{t}^{23}
$$
$$
+24214587194571650834572683444012\,{t}^{22}-
14324708604336971207868317957868\,{t}^{21}
$$
$$
+6859440538554030239641036025103\,{t}^{20}-
2654944041424536277948746010303\,{t}^{19}
$$
$$
+828959040281722890327985220255\,{t}^{18}-
208258356862493902206466194607\,{t}^{17}
$$
$$
+41958615314622858303912597215
\,{t}^{16}-6750900711491569851736413311\,{t}^{15}
$$
$$
+862879164715733847737203343\,{t}^{14}-87057778313447181201990522\,{t}^
{13}
$$
$$
+6880098547452856483997402\,{t}^{12}-422023241997789381263002\,{t}
^{11}
$$
$$
+19874722637854592209338\,{t}^{10}-709234610141846974874\,{t}^{9}
+18871718211768417242\,{t}^{8}
$$
$$
-366971376492201338\,{t}^{7}+5083730101530753\,{t}^{6}-48551559344145\,{t}^{5}+305993127089\,{t}^{4}
$$
$$
-1195782497\,{t}^{3}+2630641\,{t}^{2}-2769\,t+1.
$$

\bigskip

Note that, surprisingly, the degree of the denominator of $F_7(t)$ 
is $48$ rather than the expected $64$ since the first six generating functions' 
denominator have degree $2^{k-1}$, $1 \leq k \leq 6$. 
With a larger computer, one should be able to compute $F_k$ for larger $k$, using this experimental approach. 

Generally, for an arbitrary graph $G$, we consider the number of spanning trees in $G \times P_n$. With the same methodology, 
a list of data can be obtained empirically with which  a generating function follows.

\section*{Joint Resistance}

The original motivation for the Matrix Tree Theorem, first discovered by Kirchhoff (of Kirchhoff's laws fame)
came from the desire to efficiently compute joint resistances in an electrical network.

Suppose one is interested in the joint resistance in an electric network in the form of a grid graph between two diagonal 
vertices $[1,1]$ and $[k,n]$. We assume that each edge has resistance $1$ Ohm.
To obtain it, all we need  is, in addition for the number of spanning trees (that's the numerator),
the number of spanning forests $SF_k(n)$ of the graph $G_k(n)$ that have exactly two components, 
each component containing exactly one of the members of the pair   $\{[1,1],[k,n]\}$ (this is the denominator). 
The joint resistance is just the ratio.

In principle, we can apply the same method to obtain the generating function $S_k$. 
Empirically, we found that the denominator of $S_k$ is always the square of the denominator of $F_k$ times another polynomial $C_k$. 
Once the denominator is known, we can find the numerator in the same way as above.
So our focus is to find $C_k$.

The procedure {\tt DenomSFKN(k,t)} in the Maple package {\tt JointConductance.txt}, calculates $C_k$. For $2 \leq k \leq 4$,
we have
$$
C_2 = t-1 ,
$$
$$
C_3 = t^4-8t^3+17t^2-8t+1 ,
$$
$$
C_4 = t^{12}-46t^{11}+770t^{10}-6062t^9+24579t^8-55388t^7+72324t^6-55388t^5+24579t^4
$$
$$
-6062t^3+770t^2-46t+1 .
$$

{\bf Remark} 
By looking at the output of our Maple package, we conjectured that $R(k,n)$, the resistance between vertex $[1,1]$ and vertex $[k,n]$ in 
the $k \times n$ grid graph, $G_k(n)$, where each edge is a resistor of $1$ Ohm, 
is asymptotically $n/k$, for any fixed $k$, as $n \rightarrow \infty$. We proved it rigorously for $k \leq 6$, and
we wondered whether there is a human-generated ``electric proof''. Naturally we emailed Peter Doyle, the co-author of the
delightful masterpiece [DS], who quickly came up with the following argument.

{\it 
Making the horizontal resistors into almost 
resistance-less gold wires gives the lower bound $R(k,n) \geq (n-1)/k$ since it is a parallel circuit of $k$ resistors of $n-1$ Ohms. 
For an upper bound of the same order, put 1 Ampere in at [1,1] and out at $[k,n]$, routing $1/k$ Ampere up each of the $k$ verticals. 
The energy dissipation is $k(n-1)/k^2+C(k) = (n-1)/k+C(k)$,
where the constant $C(k)$ is the energy dissipated along the top and bottom resistors. Specifically, $C(k) = 2(1-1/k)^2 + (1-2/k)^2 + \dots + (1/k)^2)$. So
$(n-1)/k \leq R(k,n) \leq (n-1)/k + C(k)$.}

We thank Peter Doyle for his kind permission to reproduce this  {\it electrifying} argument.

\section*{The statistic of the number of vertical edges in spanning trees of grid graphs}

Often in enumerative combinatorics,  the class of interest has natural `statistics', like height, weight, and IQ for humans.
Recall that the  {\it naive counting} is
$$
|A| \, := \, \sum_{a \in A} 1 ,
$$
getting a {\bf number}. Define:
$$
|A|_x \, := \,\sum_{a \in A} x^{f(a)} ,
$$
where $f:=A \rightarrow \mathbb{Z}$ is the statistic of interest. To go  from the weighted enumeration (a certain Laurent polynomial)
to straight enumeration, one simply plugs-in $x=1$, i.e. $|A|_1 = |A|$.

The {\it scaled} random variable is defined as follows. Let $E(f)$ and $Var(f)$ be the {\it expectation} and
{\it variance}, respectively, of the statistic $f$ defined on $A$, 
and define the {\it scaled} random variable, for $a\in A$, by
$$
X (a):= \frac{f(a)- E(f)}{\sqrt{Var(f)}} .
$$
In this section, we are interested in the statistic `number of vertical  edges',
defined on spanning trees of grid graphs. 
For  given $k$ and $n$, let, as above, $G_k(n)$
denote the $k \times n$ grid-graph.
Let $\mathcal{F}_{k,n}$ be its set of spanning trees.
If the weight is 1, then $\sum_{f \in \mathcal{F}_{k,n}} 1 =|\mathcal{F}_{k,n}|$ is the naive counting. Now let's define a natural statistic 

$ver(T)$ = the number of vertical edges in the spanning tree $T$,

and the weight $w(T) =  v^{ver(T)}$, then the weighted counting follows:
$$
Ver_{k,n} (v) = \sum_{T \in \mathcal{F}_{k,n}} w(T)
$$
where $\mathcal{F}_{k,n}$ is the set of spanning trees of $G_k(n)$. 

We define the bivariate generating function
$$
g_{k}(v,t) = \sum_{n=0}^{\infty} Ver_{k,n} t^n.
$$
More generally, with our Maple package {\tt GFMatrix.txt}, and procedure {\tt VerGF}, 
we are able to obtain the bivariate generating function for an arbitrary graph of the form $G \times P_n$. 
The procedure {\tt VerGF} takes inputs $G$ (an arbitrary graph), $N$ (an integer determining how many data we use to 
find the recurrence relation) and two symbols $v$ and $t$.

The main tool for computing {\tt VerGF} is still the 
Matrix Tree Theorem and {\tt GuessRec}. 
But we need to modify  the Laplacian matrix for the graph. 
Instead of letting $a_{ij}=-1$ for $i \neq j$ and $\{i,j\} \in E(G \times P_n)$, 
we should consider whether the edge $\{i,j\}$ is a vertical edge. 
If so, we let $a_{i,j}=-v, a_{j,i}=-v$. The diagonal elements which are $(-1) \times$ 
(the sum of the rest entries on the same row) should change accordingly. 
The following theorems are for grid graphs when $2 \leq k \leq 4$ while $k=1$ is a trivial case because there are no vertical edges.

\bigskip

{\bf Theorem 8} The bivariate generating function for the weighted counting according to the number of vertical edges of spanning trees in $G_2(n)$ is:
$$
g_2(v,t) = \frac {vt}{1- \left( 2\,v+2 \right) t+{t}^{2}}  .
$$

\bigskip

{\bf Theorem 9} The bivariate generating function for the weighted counting 
 according to the number of vertical edges  vertical edges of spanning trees in $G_3(n)$ is:
$$
g_3(v,t) = \frac {-{t}^{3}{v}^{2}+{v}^{2}t}{1- \left( 3\,{v}^{2}+8\,v+4 \right) t- \left( -10\,{v}^{2}-16\,v-6 \right) {t}^{2}- \left( 3\,{v}^{2}+8\,v+
4 \right) {t}^{3}+{t}^{4}} .
$$

\bigskip

{\bf Theorem 10} The bivariate generating function for the weighted counting
 according to the number of vertical edges of spanning trees in $G_4(n)$ is:
$$
g_4(v,t) = \frac{numer(g_4)}{denom(g_4)}
$$
where
$$
numer(g_4) = {v}^{3}t+ \left( -16\,{v}^{5}-24\,{v}^{4}-9\,{v}^{3} \right) {t}^{3}+
 \left( 8\,{v}^{6}+40\,{v}^{5}+48\,{v}^{4}+16\,{v}^{3} \right) {t}^{4}
$$
$$
+ \left( -16\,{v}^{5}-24\,{v}^{4}-9\,{v}^{3} \right) {t}^{5}+{v}^{3}{t
}^{7}
$$
and
$$
denom(g_4) = 1- \left( 4\,{v}^{3}+20\,{v}^{2}+24\,v+8 \right) t- \left( -52\,{v}^{4
}-192\,{v}^{3}-256\,{v}^{2}-144\,v-28 \right) {t}^{2}
$$
$$
- \left( 64\,{v}^
{5}+416\,{v}^{4}+892\,{v}^{3}+844\,{v}^{2}+360\,v+56 \right) {t}^{3}
$$
$$
- \left( -16\,{v}^{6}-160\,{v}^{5}-744\,{v}^{4}-1408\,{v}^{3}-1216\,{v}
^{2}-480\,v-70 \right) {t}^{4}
$$
$$
- \left( 64\,{v}^{5}+416\,{v}^{4}+892\,{
v}^{3}+844\,{v}^{2}+360\,v+56 \right) {t}^{5}- \left( -52\,{v}^{4}-192
\,{v}^{3}-256\,{v}^{2}-144\,v-28 \right) {t}^{6}
$$
$$
- \left( 4\,{v}^{3}+20
\,{v}^{2}+24\,v+8 \right) {t}^{7}+{t}^{8}  .
$$
With the Maple package {\tt BiVariateMoms.txt} and its {\tt Story} procedure from

\bigskip

{\tt http://sites.math.rutgers.edu/\~{}zeilberg/tokhniot/BiVariateMoms.txt},

\bigskip

the expectation, variance and higher moments can be easily analyzed. We calculated up to the 4th moment for $G_2(n)$. 
For $k=3,4$, you can find the output files from

{\tt http://sites.math.rutgers.edu/\~{}yao/OutputStatisticVerticalk=3.txt}

{\tt http://sites.math.rutgers.edu/\~{}yao/OutputStatisticVerticalk=4.txt}

\bigskip

{\bf Theorem 11} The moments of the statistic: the number of vertical edges in the spanning trees of $G_2(n)$ are as follows:

Let $b$ be the largest positive root of the polynomial equation
$$
b^2-4b+1 = 0
$$
whose floating-point approximation is 3.732050808, then the size of the $n$-th family (i.e. straight enumeration) is very close to
$$
{\frac {{b}^{n+1}}{-2+4\,b}} .
$$
The average of the statistics is, asymptotically 
$$
\frac{1}{3}+\frac{1}{3}\,{\frac { \left( -1+2\,b \right) n}{b}}  .
$$
The variance of the statistics is, asymptotically 
$$
-\frac{1}{9}+\frac{1}{9}\,{\frac { \left( 7\,b-2 \right) n}{-1+4\,b}} .
$$
The skewness of the statistics is, asymptotically 
$$
{\frac {780\,b-209}{ \left( 4053\,b-1086 \right) {n}^{3}+ \left( -7020
\,b+1881 \right) {n}^{2}+ \left( 4053\,b-1086 \right) n-780\,b+209}}.
$$
The kurtosis of the statistics is, asymptotically 
$$
3\,{\frac { \left( 32592\,b-8733 \right) {n}^{2}+ \left( -56451\,b+
15126 \right) n+21728\,b-5822}{ \left( 32592\,b-8733 \right) {n}^{2}+
 \left( -37634\,b+10084 \right) n+10864\,b-2911}} .
$$

\section*{Application of the C-finite Ansatz to computing generating functions of determinants (and permanents) of  almost-diagonal  Toeplitz matrices}

So far, we have seen applications of the $C$-finite ansatz methodology 
for automatically computing generating functions for
enumerating spanning trees/forests for certain infinite families of graphs. 

The second case study is completely different, and in a sense more general,
since the former framework may be subsumed in this new context.

\bigskip

{\bf Definition} Diagonal matrices $A$ are square matrices in which the entries outside the main diagonal are $0$, i.e. $a_{ij} = 0$ if $i \neq j$.

\bigskip

{\bf Definition}  An almost-diagonal  Toeplitz matrix $A$ is a square matrices 
in which $a_{i,j} = 0$ if $j-i \geq k_1$ or $i-j \geq k_2$ for some fixed positive integers $k_1, k_2$ 
and $\forall i_1, j_1, i_2, j_2$, if $i_1-j_1 = i_2-j_2$, then $a_{i_1 j_1} = a_{i_2 j_2}$. 

\bigskip

For simplicity, we use the notation $L=$[n, [the first $k_1$ entries in the first row], [the first $k_2$ entries in the first column]] to denote 
the $n \times n$ matrix with these specifications.
Note that this notation already contains all information we need to reconstruct this matrix. 
For example, [6, [1,2,3], [1,4]] is the matrix
$$
\begin{bmatrix}
 1 & 2 & 3 & 0 & 0 & 0 \\
 4 & 1 & 2 & 3 & 0 & 0 \\
 0 & 4 & 1 & 2 & 3 & 0 \\
 0 & 0 & 4 & 1 & 2 & 3 \\
 0 & 0 & 0 & 4 & 1 & 2 \\
 0 & 0 & 0 & 0 & 4 & 1
\end{bmatrix} .
$$

The following is the Maple procedure {\tt DiagMatrixL}  (in our Maple package {\tt GFMatrix.txt}),
which inputs such a list $L$ and outputs the corresponding matrix. 

\bigskip

{\obeylines
{\tt
DiagMatrixL:=proc(L) local n, r1, c1,p,q,S,M,i:
n:=L[1]:
r1:=L[2]:
c1:=L[3]:
p:=nops(r1)-1:
q:=nops(c1)-1:
if r1[1] <> c1[1] then
\quad   return fail:
fi:
S:=[0\$(n-1-q), seq(c1[q-i+1],i=0..q-1), op(r1), 0\$(n-1-p)]:
M:=[0\$n]:
for i from 1 to n do
\quad   M[i]:=[op(max(0,n-1-q)+q+2-i..max(0,n-1-q)+q+1+n-i,S)]:
od:
return M:
end:
}
}

\bigskip

For this matrix, $k_1=3$ and $k_2=2$. 
Let $k_1, k_2$ be fixed and $M_1, M_2$ be two lists of numbers or symbols of length $k_1$ and $k_2$ respectively, 
$A_k$ is the almost-diagonal  Toeplitz matrix represented by the list $L_k = [k, M_1, M_2]$. Note that the first elements in the lists $M_1$ and $M_2$ must be identical. 

Having fixed two lists $M_1$ of length $k_1$ and $M_2$ of length $k_2$, (where $M_1[1]=M_2[1]$), it is of interest
to derive {\it automatically}, the generating function (that is always a rational function for reasons that will
soon become clear), $\sum_{k=0}^{\infty} a_k \, t^k$, where
 $a_k$ denotes the determinant of the $k \times k$ almost-diagonal  Toeplitz matrix whose first row starts with $M_1$, and first column
starts with $M_2$. Analogously, it is also of interest to do the analogous problem when the determinant is replaced by the permanent.

Here is the Maple procedure {\tt GFfamilyDet} which takes inputs (i) $A$: a name of a Maple procedure that inputs an integer $n$ 
and outputs an $n \times n$ matrix according to some rule, e.g., the almost-diagonal  Toeplitz matrices, 
(ii) a variable name $t$, (iii) two integers $m$ and $n$ which are the lower and upper bounds 
of the sequence of determinants we consider. It outputs a rational function in $t$, say $R(t)$, which is the generating function of the sequence.
\vspace{1mm}\noindent
{\obeylines
{\tt
GFfamilyDet:=proc(A,t,m,n) local i,rec,GF,B,gu,Denom,L,Numer:
L:=[seq(det(A(i)),i=1..n)]:
rec:=GuessRec([op(m..n,L)])[2]:
gu:=solve(B-1-add(t**i*rec[i]*B,i=1..nops(rec)), {B}):
Denom:=denom(subs(gu,B)):
Numer:=Denom*(1+add(L[i]*t**i, i=1..n)):
Numer:=add(coeff(Numer,t,i)*t**i, i=0..degree(Denom,t)):
Numer/Denom:
end: 
}
}
\vspace{1mm}\noindent
Similarly we have procedure {\tt GFfamilyPer} for the permanent. Let's look at an example. 
The following is a sample procedure which considers the family of 
almost diagonal  Toeplitz matrices which the first row $[2,3]$ and the first column $[2,4,5]$. 
\vspace{1mm}\noindent
{\obeylines
{\tt
SampleB:=proc(n) local L,M:
L:=[n, [2,3], [2,4,5]]:
M:=DiagMatrixL(L):
end:
}
}
Then {\tt GFfamilyDet(SampleB, t, 10, 50)} will return the generating function 
$$
-\frac{1}{ 45\,{t}^{3}-12\,{t}^{2}+2\,t-1 } .
$$

It turns out, that for this problem, the more `conceptual' approach of setting up a transfer matrix
also works well. But don't worry, the computer can do the `research' all by itself, with only
a minimum amount of human pre-processing.

We will now describe this more conceptual approach, that may be called {\it symbolic dynamical programming}, where
the computer sets up, {\it automatically}, a finite-state scheme, by {\it dynamically} discovering the
set of states, and automatically figures out the transfer matrix.

\section*{The Transfer Matrix  method for almost-diagonal  Toeplitz matrices}

Recall from Linear Algebra 101, the

{\bf Cofactor Expansion} Let $|A|$ denote the determinant of an $n \times n$ matrix $A$, then
$$
|A| = \sum_{j=1}^{n} (-1)^{i+j} a_{ij} M_{ij}, \quad \forall i \in [n],
$$
where $M_{ij}$ is the $(i,j)$-minor.

We'd like to consider the Cofactor Expansion for almost-diagonal  Toeplitz matrices along the first row. 
For simplicity, we assume while $a_{i,j} = 0$ if $j-i \geq k_1$ or $i-j \geq k_2$ 
for some fixed positive integers $k_1, k_2$, and if $-k_2 < j_1-i_1 < j_2-i_2 < k_1$, then $a_{i_1 j_1} \neq a_{i_2 j_2}$. 
Under this assumption, for any minors we obtain through recursive Cofactor Expansion along the first row, 
the dimension, the first row and the first column should provide enough information to reconstruct the matrix. 

For an almost-diagonal  Toeplitz matrix represented by $L=$[Dimension, [the first $k_1$ entries in the first row], 
[the first $k_2$ entries in the first column]], any minor can be represented by 
[Dimension, [entries in the first row up to the last nonzero entry], [entries in the first column up to the last nonzero entry]].

Our goal in this section is the same as the last one, to get a generating function for the determinant or permanent of almost-diagonal  Toeplitz matrices $A_k$ with dimension $k$. Once we have 
those almost-diagonal  Toeplitz matrices, the first step is to do a one-step expansion as follows:
\vspace{1mm}\noindent
{\obeylines
{\tt
ExpandMatrixL:=proc(L,L1) 
local n,R,C,dim,R1,C1,i,r,S,candidate,newrow,newcol,gu,mu,temp,p,q,j:
n:=L[1]:
R:=L[2]:
C:=L[3]:
p:=nops(R)-1:
q:=nops(C)-1:
dim:=L1[1]:
R1:=L1[2]:
C1:=L1[3]:

if R1=[] or C1=[] then
\quad  return {}:
elif R[1]<>C[1] or R1[1]<>C1[1] or dim>n then 
\quad  return fail:
else

  S:=\{\}:

  gu:=[0\$(n-1-q), seq(C[q-i+1],i=0..q-1), op(R), 0\$(n-1-p)]:
  candidate:=[0\$nops(R1),R1[-1]]:
  for i from 1 to nops(R1) do
    \quad mu:=R1[i]:
    for j from n-q to nops(gu) do
     \quad \quad if gu[j]=mu then
       \quad\quad\quad  candidate[i]:=gu[j-1]:
\quad \quad      fi:
\quad    od:
  od:
      
  for i from n-q to nops(gu) do
  \quad  if gu[i] = R1[2] then
     \quad \quad temp:=i:
      \quad \quad break:
  \quad    fi:
  od:

  for i from 1 to nops(R1) do
\quad    if i = 1 then
\quad\quad       mu:=[R1[i]*(-1)**(i+1), [dim-1,[op(i+1..nops(candidate), candidate)], [seq(gu[temp-i],i=1..temp-n+q)]]]:
\quad\quad      S:=S union {mu}:
    \quad else
      \quad \quad mu:=[R1[i]*(-1)**(i+1), [dim-1, [op(1..i-1, candidate), op(i+1..nops(candidate), candidate)], [op(2..nops(C1), C1)]]]:
\quad\quad      S:=S union {mu}:
\quad    fi:
  od:

\quad  return S:

fi:

end:
}
}
\vspace{1mm}\noindent
The {\tt ExpandMatrixL} procedure inputs a data structure $L =$ [Dimension, first\_row=[ ], first\_col=[ ]] 
as the matrix we start and the other data structure $L1$ as the current minor we have, 
expands $L1$ along its first row and outputs a list of [multiplicity, data structure]. 

We would like to generate all the "children" of an almost-diagonal  Toeplitz matrix regardless of the dimension, i.e., 
two lists $L$ represent the same child as long as their first\_rows and first\_columns are the same, respectively. 
The set of "children" is the scheme of the almost diagonal  Toeplitz matrices in this case.

The following is the Maple procedure {\tt ChildrenMatrixL} which inputs a data structure $L$ and 
outputs the set of its "children" under Cofactor Expansion along the first row:

\bigskip

{\obeylines
{\tt
ChildrenMatrixL:=proc(L) local S,t,T,dim,U,u,s:
dim:=L[1]:
S:=\{[op(2..3,L)]\}:
T:=\{seq([op(2..3,t[2])],t in ExpandMatrixL(L,L))\}:
while T minus S <> \{\} do
\quad  U:=T minus S:
\quad  S:=S union T:
\quad   T:=\{\}:
\quad   for u in U do
\quad    T:=T union \{seq([op(2..3,t[2])],t in ExpandMatrixL(L,[dim,op(u)]))\}:
  od:
od:
for s in S do
\quad  if s[1]=[] or s[2]=[] then
\quad\quad    S:=S minus \{s\}:
\quad  fi:
od:
S:
end:
}
}

\bigskip

After we have the scheme $S$, by the Cofactor Expansion of any element in the scheme, a system of algebraic equations follows. 
For children in $S$, it's convenient to let the almost-diagonal  Toeplitz matrix be the first one $C_1$ and for 
the rest, any arbitrary ordering will do. For example, if after Cofactor Expansion for $C_1$, $c_2$ "copies" of 
$C_2$ and $c_3$ "copies" of $C_3$ are obtained, then the equation will be 
$$
C_1 = 1+ c_2 t C_2 + c_3 t C_3  .
$$
However, if the above equation is for $C_i, i \neq 1$, i.e. $C_i$ is not the almost-diagonal  Toeplitz matrix itself, then the equation will be slightly different:
$$
C_i = c_2 t C_2 + c_3 t C_3 .
$$
Here $t$ is a symbol as we assume the generating function is a rational function of $t$.

Here is the Maple code that implements how we get the generating function for the determinant of a family of 
almost-diagonal  Toeplitz matrices by solving a system of algebraic equations:

\bigskip

{\obeylines
{\tt
GFMatrixL:=proc(L,t) local S,dim,var,eq,n,A,i,result,gu,mu:
dim:=L[1]:
S:=ChildrenMatrixL(L):
S:=[[op(2..3,L)], op(S minus \{[op(2..3,L)]\})]:
n:=nops(S):
var:=\{seq(A[i],i=1..n)\}:
eq:=\{\}:
for i from 1 to 1 do
\quad   result:=ExpandMatrixL(L,[dim,op(S[i])]):
\quad  for gu in result do
\quad \quad  if gu[2][2]=[] or gu[2][3]=[] then
\quad \quad \quad   result:=result minus \{gu\}:
\quad \quad fi:
\quad  od:
\quad  eq:=eq union \{A[i] - 1 - add(gu[1]*t*A[CountRank(S, [op(2..3, gu[2])])], gu in result)\}:
od:
for i from 2 to n do
\quad  result:=ExpandMatrixL(L,[dim,op(S[i])]):
\quad  for gu in result do
\quad  if gu[2][2]=[] or gu[2][3]=[] then
\quad \quad    result:=result minus {gu}:
\quad  fi:
  od:
  eq:=eq union \{A[i] - add(gu[1]*t*A[CountRank(S, [op(2..3, gu[2])])], gu in result)\}:
od:
gu:=solve(eq, var)[1]:
subs(gu, A[1]):

end:

}
}

\bigskip

{\tt GFMatrixL([20, [2, 3], [2, 4, 5]], t)} returns
$$
- \frac{1}{ 45\,{t}^{3}-12\,{t}^{2}+2\,t-1}  .
$$
Compared to empirical approach, the `symbolic dynamical programming' method is faster and more efficient
for the moderate-size examples that we tried out.
However, as the lists will grow larger, it is likely that the former method will win out,
since with this non-guessing approach, it is equally fast to get generating functions for
determinants and permanents, and as we all know, permanents are hard.

The advantage of the present method  is that it is more appealing to humans, and does not
require any `meta-level' act of faith.
However, both methods are very versatile and are great experimental approaches for enumerative combinatorics problems. 
We hope that our readers will find other applications.

\section*{Summary}
Rather than trying to tackle each enumeration problem, one at a time, using
ad hoc human ingenuity each time, building up an intricate transfer matrix, and
only using the computer at the end as a symbolic calculator, it is a much better
use of our beloved silicon servants
(soon to become our masters!) to replace `thinking' by `meta-thinking', i.e. develop experimental mathematics methods
that can handle many different types of problems. In the two case studies discussed here,
everything was made rigorous, but if one can make semi-rigorous and even non-rigorous
discoveries, as long as they are {\it interesting}, one should not be hung up
on rigorous proofs. In other words, if you can find a rigorous justification (like in these
two case studies) that's nice, but if you can't, that's also nice!

\section*{Acknowledgment} Many thanks are due to a very careful referee that pointed out many minor, but annoying,  errors.
Also thanks to Peter Doyle for permission to include his elegant electric argument.

\section*{References}

[DS] 
Peter Doyle and Laurie Snell, {\it ``Random Walks and Electrical Networks''},
Carus Mathematical Monographs (\# 22), Math. Assn. of America, 1984.

\bigskip

[ESZ] Shalosh B. Ekhad, N. J. A. Sloane and Doron Zeilberger,
{\it Automated Proof (or Disproof) of Linear Recurrences Satisfied by Pisot Sequences}, 
The Personal Journal of Shalosh B. Ekhad and Doron Zeilberger.
Available from  \hfill\break
{\tt http://www.math.rutgers.edu/\~{}zeilberg/mamarim/mamarimhtml/pisot.html} 
\bigskip

[F] F. J. Faase, {\it On the number of specific spanning subgraphs of the graphs $g \times p_n$}, Ars Combinatorica, {\bf 49} (1998)
129-154

\bigskip

[KP] Manuel Kauers  and Peter Paule, {\it ``The Concrete Tetrahedron''}, Springer, 2011.

\bigskip

[R] Paul Raff, {\it Spanning Trees in Grid Graph}, \hfill\break
{\tt https://arxiv.org/abs/0809.2551}. 

\bigskip

[S] Richard Stanley, {\it ``Enumerative Combinatorics, Volume 1''}.
First edition: Wadsworth \& Brooks/Cole, 1986.
Second edition: Cambridge University Press, 2011.

\bigskip

[Wi] Wikipedia contributors. ``Berlekamp-Massey algorithm.'' Wikipedia, The Free Encyclopedia. Wikipedia, The Free Encyclopedia, 26 Nov. 2018. Web. 7 Jan. 201

\bigskip

[Z1] Doron Zeilberger, {\it The C-finite Ansatz},  Ramanujan Journal {\bf 31} (2013), 23-32. 
Available from  \hfill\break
{\tt http://www.math.rutgers.edu/\~{}zeilberg/mamarim/mamarimhtml/cfinite.html} 

\bigskip

[Z2] Doron Zeilberger, {\it Why the Cautionary Tales Supplied by Richard Guy's Strong Law of Small Numbers Should not be Overstated},
The Personal Journal of Shalosh B. Ekhad and Doron Zeilberger.
Available from  \hfill\break
{\tt http://www.math.rutgers.edu/\~{}zeilberg/mamarim/mamarimhtml/small.html}

\bigskip

[Z3] Doron Zeilberger,
{\it An Enquiry Concerning Human (and Computer!) [Mathematical] Understanding},
in: C.S. Calude, ed., ``Randomness \& Complexity, from Leibniz to Chaitin'' World Scientific, Singapore, 2007, pp. 383-410.
Available from  \hfill\break
{\tt http://www.math.rutgers.edu/\~{}zeilberg/mamarim/mamarimhtml/enquiry.html}

\bigskip

First Version: Dec. 17, 2018. This version: Jan. 8, 2019.

\end{document}